\documentclass[11pt, leqno]{article}
\usepackage{latexsym,amssymb,amsthm}
\usepackage{amsmath}
\usepackage{color}

\usepackage{multirow,eepic}

\makeatletter
\def\senbun#1(#2)#3({\@senbun(#2)(}
\def\@senbun(#1,#2)(#3,#4){%
   \@tempdima#1\p@ \advance\@tempdima#3\p@
   \divide\@tempdima\tw@
   \@tempdimb#2\p@ \advance\@tempdimb#4\p@
   \divide\@tempdimb\tw@
   \edef\@senbun@temp{\noexpand\qbezier(#1,#2)%
      (\strip@pt\@tempdima,\strip@pt\@tempdimb)(#3,#4)}%
   \@senbun@temp}
\makeatother

\newtheorem{definition}{Definition}[section]

\newtheorem{lemma}[definition]{Lemma}

\newtheorem{theorem}[definition]{Theorem}

\newtheorem{proposition}[definition]{Proposition}

\newtheorem{remark}[definition]{Remark}

\makeatletter
\long\def\unmarkedfootnote#1{{\long\def\@makefntext##1{##1}\footnotetext{#1}}}
\makeatother

\renewcommand{\qed}{\qquad\kern1pt   
   \vbox{\hrule height 0.6pt      
         \hbox{\vrule width 0.6pt 
               \vbox{\vskip 6pt}  
               \hskip 3pt
              \vrule width 1.3pt} 
         \hrule depth 1.3pt}     
   \kern1pt}

\numberwithin{equation}{section}

\begin{document}

\title{
Canceling  effects in  higher-order  Hardy-Sobolev inequalities }
\author{
Andrea Cianchi\\
Dipartimento di Matematica e Informatica, Universit\`{a} di Firenze\\
Viale Morgagni 67/A, Firenze 50134, Italy
\vspace{5pt}\\
\\
Norisuke Ioku\\
Graduate School of Science and Engineering, Ehime University\\
Matsuyama, Ehime 790-8577, Japan\vspace{5pt}\\
}
\date{}
\maketitle

\begin{abstract}
%
A classical   first-order  Hardy-Sobolev inequality in Euclidean
 domains, involving weighted
  norms   depending on powers of the distance function from
 their boundary, is known to hold for every, but one, value of the
 power.
 We show that, by contrast, the  missing power is  admissible in a
 suitable counterpart for higher-order Sobolev norms. 
 Our result
 complements and extends  contributions by Castro and Wang \cite{CW}, and Castro, D\'avila and Wang \cite{CDW1, CDW2},  where a surprising
 canceling phenomenon underling the relevant inequalities was
 discovered in the special case of functions with derivatives in
 $L^1$. 
\end{abstract}

\unmarkedfootnote {\par\noindent {\it Mathematics Subject
Classification: 46E35, 46E30.}
}

\section{Introduction and main results}

%


Weighted Sobolev inequalities, namely Sobolev inequalities for norms
on open sets $\Omega$ of $\mathbb R^n$ equipped with measures having
densities -- the weights -- with respect to the Lebesgue measure,
have been extensively investigated, mainly in connection with  the
theory of degenerate partial differential equations.
The literature in this area is very
rich.
%
Let us just mention that various characterizations of the weights
supporting the relevant inequalities are available, such as those
depending on global integrability properties of the weights
\cite{MS}, on the growth of their integrals on balls \cite{Ad1,
Ad2}, on their membership to Muckenhoupt $A_p$-classes \cite{HKM,
OK}, on associated capacities \cite{mazya}, on their rearrangements
\cite{CEG}.
\par The most popular and widely  exploited inequalities of this kind are presumably the so called
Hardy-Sobolev inequalities, whose weights are just powers of the
distance function from the boundary  $\partial \Omega$.
This
function
%
is known to inherit regularity properties of $\partial \Omega$ in a
sufficiently narrow neighborhood of the latter. Let us call $d :
\Omega \to (0, \infty)$ a function which agrees with the distance
function in such neighborhood of $\partial \Omega$, and enjoys the
same regularity properties, but in the whole of $\Omega$.
\\
Given $p \in [1, \infty]$ and $\alpha \in \mathbb R$, we denote by
$L^p(\Omega , d^\alpha)$ the weighted Lebesgue space equipped with
the norm defined as
$$\|u\|_{L^p(\Omega , d^\alpha)} = \bigg(\int_\Omega |u(x)|^p
d(x)^\alpha \, dx\bigg)^{\frac 1p}$$ for a measurable function $u$
in $\Omega$. Moreover, if $m \in \mathbb N$, the notation
$W^{m,p}(\Omega , d^\alpha)$ is adopted for the associated Sobolev
space of $m$-times weakly differentiable functions $u$ in $\Omega$
endowed with the norm
$$\|u\|_{W^{m,p}(\Omega , d^\alpha)}= \sum _{j=0}^m \|\nabla ^j u\|_{L^p(\Omega ,
d^\alpha)},$$ where $\nabla ^j u$ stands for the vector of all
derivatives of $u$ of order $j$. We also simply denote $\nabla ^1 u$
by $\nabla u$; also, $\nabla ^0 u$ stands for $u$. The notation
$W^{m,p}_0(\Omega , d^\alpha)$ is devoted to the closure of
$C^\infty _0(\Omega)$ in $W^{m,p}(\Omega , d^\alpha)$.
\par A classical  Hardy-Sobolev  inequality  asserts that if
$\Omega$ is a bounded Lipschitz domain, and
 $$\alpha \neq p-1,$$
then there exists a constant $C$ such that
\begin{equation}\label{classical}
\left\|\frac u{d} \right\|_{L^p(\Omega, d^{\alpha})} \leq C \|u\|_{W
^{1,p}(\Omega , d^{\alpha})}
\end{equation}
for every $u \in W_0^{1,p}(\Omega , d^{\alpha})$ \cite[Theorem~8.4]{K}. On the other hand, inequality \eqref{classical} fails for
the critical value $\alpha = p-1$.
\par
The main purpose of the present paper is to show that, this
notwithstanding, suitable higher-order versions of inequality
\eqref{classical}, which cannot just be obtained  from
\eqref{classical} via iteration, do hold  even when $\alpha = p-1$.
\\ A prototypal  second-order
inequality may help to grasp the spirit of our results. Assume that
$\Omega$ has a smooth boundary, so that $d$ is also smooth in a
neighborhood of $\partial \Omega$.
%
%
 Let $u \in W_0^{2,p}(\Omega , d^{p-1})$. A standard
 property of the distance function ensures that
  $|\nabla d|=1$  in a
neighborhood of $\partial \Omega$, whence
\begin{equation}\label{cancel}
\left|\nabla \Big(\frac ud\Big)\right| = \left|\frac{\nabla u} d - u
\frac{\nabla d}{d^2}\right| \leq  \frac{|\nabla u|}{d} + \frac
{|u|}{d^2}
\end{equation}
a.e. in the same neighborhood. Inequality \eqref{classical} cannot
be exploited to infer that the functions $\frac{|\nabla u|}{d}$ and
$\frac {|u|}{d^2}$ belong to $L^p(\Omega , d^{p-1})$. In fact,
membership to $L^p(\Omega , d^{p-1})$ of neither of these functions
is guaranteed under
the sole assumption that $u \in W_0^{2,p}(\Omega , d^{p-1})$
(this can be verified, for instance, by taking $\Omega = (0,1)$, and
considering functions $u(x)$ decaying like $x\log ^{-\alpha} (\tfrac
1x)$ as $x \to 0^+$, with $\alpha \in (0, \tfrac 1p]$).
 Nevertheless, we show that the
inequality
\begin{equation}\label{second}
\left\|\frac u{d} \right\|_{W^{1,p} (\Omega, d^{p-1})} \leq C
\|u\|_{W ^{2,p}(\Omega , d^{p-1})},
\end{equation}
 holds for some
 constant $C$, and for every $u \in W_0^{2,p}(\Omega , d^{p-1})$.
 This is possible thanks to a canceling effect
which allows the leftmost side of \eqref{cancel} to have stronger
integrability properties than each  addend on its rightmost side.
%
%
%
\\
In the case when $p=1$, such a striking phenomenon has been elucidated in  remarkable
 contributions,  by which ours is inspired, of Castro and Wang \cite{CW}, for $n=1$, and of Castro,
 D\'avila and Wang \cite{CDW1, CDW2}, for any $n\geq 1$.
 In this case,
 non-weighted Lebesgue and Sobolev norms appear in
 \eqref{classical} and \eqref{second}, and in their higher-order
 counterparts   from \cite{CDW2}. 
\par The arbitrary-order
version of inequality \eqref{second} to be established asserts that,
if $1\le p <\infty$,  and $k,m\in \mathbb{N}$, with $m\ge 2$,
 and $1\le k\le m-1$, then there exists a
 constant
 $C$ such that
\begin{equation}\label{eq:1.3p}
\left\|\frac u{d^{m-k}} \right\|_{W^{k,p} (\Omega, d^{p-1})} \leq C
\|u\|_{W ^{m,p}(\Omega , d^{p-1})}
\end{equation}
for every $u \in W_0^{m,p}(\Omega , d^{p-1})$.
\par Inequality
\eqref{eq:1.3p} is in turn a special instance of our most general
result, stated in the following theorem, where Sobolev type spaces
associated with different Lebesgue norms and distance weights are
allowed on the two sides of the relevant inequality.


\begin{theorem}\label{proposition:1.3}
Let $\Omega$ be a bounded  open set with a smooth boundary in
$\mathbb R^n$, $n \geq 1$,  and let  $k,m\in \mathbb{N}$,  $m\ge 2$,
 and $1\le k\le m-1$.
 Assume  that
 $1\le p \le q <\infty$, and
 \begin{equation}\label{pq}
\frac 1q \geq \frac {n-p(m-k)}{np}.
\end{equation}
Let
\begin{equation}\label{r}
r \geq  \frac {q (n-1) -p(n -q)}p\,.
\end{equation}
%
%
  Then, there exists a
constant
 $C$ such that
 \begin{equation}\label{eq:1.3}
\left\|\frac u{d^{m-k}}\right\|_{W^{k,q}(\Omega, d^r)} \leq C
\|u\|_{W^{m,p}(\Omega , d^{p-1})}
\end{equation}
for every $u \in W_0^{m,p}(\Omega , d^{p-1})$.
\end{theorem}

\begin{remark}\label{rem1}
{\rm Conditions  \eqref{pq} and \eqref{r} in Theorem
\ref{proposition:1.3} are sharp, as shown in Propositions
\ref{propo2} and \ref{propo1}, Section \ref{bounded}. The assumption
that $k \geq 1$ is also sharp, since inequality \eqref{eq:1.3}
breaks down for $k=0$, as pointed out in the discussion  above.
%
%
}
\end{remark}

\begin{remark}\label{remp=1}
{\rm As already mentioned, the case when $p=q=1$ and $r=0$ in
Theorem \ref{proposition:1.3} is the object of \cite{CDW2}. Its
one-dimensional version was earlier proved in \cite{CW}, where a
higher-order inequality for $p=q>1$ is also established. However,
that higher-order inequality does not correspond to the critical
missing case of \eqref{classical}, and can be derived through  a
repeated application of the latter.}
\end{remark}

Our proof of Theorem
 \ref{proposition:1.3} combines a flattening argument for $\partial
 \Omega$, which was introduced in \cite{CDW2} and involves highly original tricks, with Whitney type decomposition techniques
 exploited in \cite{H1, H2}, and with one-dimensional Hardy type
 inequalities. In comparison with the proofs of \cite{CDW2}, additional
 difficulties arise, due to the presence of weights and of possibly
 different norms on the two sides of the inequalities under
 consideration. In particular, an iterative argument relying upon the
 use of a fundamental second-order inequality as in \cite{CDW2} is not
 possible.
 This calls for a different, direct proof, which requires a
 careful combinatorial analysis of the mutual canceling of partial
 derivatives of trial functions. 


\section{Inequalities in the half space
$\mathbb{R}_+^n$
}\label{halfspace}

This section is devoted to a Hardy-Sobolev inequality in the
half-space, contained in Theorem \ref{proposition:1.1} below. This
is a key step, of independent interest, towards the proof of Theorem
\ref{proposition:1.3}.

\begin{theorem}\label{proposition:1.1}
  Let $k,m\in \mathbb{N}$, $m\ge 2$ and $1\le k\le m-1$.
 Assume that
 $1\le p \le q <\infty$,
 $\frac1p-\frac1q= \frac{\beta-\alpha}{n}$,
 $\alpha<k+
\frac{p-1}{p}
 $,
 and $\alpha \le \beta \le \alpha+(m-k)$.
 Then, there exists a
  constant
 $C$ such that
  \begin{equation}\label{eq:1.1}
\left(
 \int_{\mathbb{R}_+^n}x_n^{\beta q}
  \left|
   \nabla ^k \left(\frac{u(x)}{x_n^{m-k}}\right)
  \right|^qdx
\right)^{\frac1q} \le C \left(
 \int_{\mathbb{R}_+^n}x_n^{\alpha p}
  \left|
   \nabla ^m u
  \right|^pdx
\right)^{\frac{1}{p}}
\end{equation}
for every $u\in C_0^{\infty}(\mathbb{R}^{n}_+)$.
 \end{theorem}

The proof of Theorem \ref{proposition:1.1} rests upon several
lemmas. The first one is a one-dimensional version of its, in the
special case when $p=q$ and $\alpha = \beta$.

\begin{lemma}\label{lemma:2.1}
 Let $k,m\in \mathbb{N}$,  $m\ge 2$ and $1\le k\le m-1$.
Assume that $1\le p<\infty$ and $\alpha<k+\frac{p-1}{p}$. Then there
exists a
constant $C$ such that
\begin{equation}\label{eq:2.1}
 \left(
 \int_0^{\infty}
  r^{\alpha p}
   \left|
    \frac{d^k}{dr^k}\left(\frac{f(r)}{r^{m-k}}\right)
   \right|^p
   dr
 \right)^{\frac{1}{p}}
\le
C
 \left(
 \int_0^{\infty}
  r^{\alpha p}
   \left|
     \frac{d^mf}{dr^m}   (r)
   \right|^p
   dr
 \right)^{\frac{1}{p}}
\end{equation}
for every $f\in C_0^{\infty}(0,\infty)$.
\end{lemma}

\noindent {\bf Proof}. An application of Taylor's formula, with
remainder term in integral form, tells us that
\begin{equation}\label{eq:2.2}
 \frac{d^k}{dr^k}
  \left(
   \frac{f(r)}{r^{m-k}}
  \right)
   =
\frac{1}{(m-k-1)!}   \int_0^r  \frac{d^m f}{ds^m}  (s)
    \left(
     1-\frac{s}{r}
    \right)^{m-k-1}
    \left(
     \frac{s}{r}
    \right)^{k-1}
   \frac{s}{r^2}ds \quad \hbox{for $r>0$,}
\end{equation}
see
\cite[Proof of Theorem
1.2]{CW}.
Hence, since $ \left(
     1-\frac{s}{r}
    \right)^{m-k-1}
\le 1$ if $0\leq s \leq r$,  one has that
\begin{equation}\label{april1}
 \left|
   \frac{d^k}{dr^k}\left(\frac{f(r)}{r^{m-k}}\right)
 \right|
 \le \frac{1}{r^{(k+1)}(m-k-1)!} \int_0^r s^k\left|
 \frac{d^m f}{ds^m}  (s)\right|ds \quad \hbox{for $r>0$.}
\end{equation}
Owing to the assumption that $\alpha<k+ \frac{p-1}{p}$, inequalities
\eqref{eq:2.2} and \eqref{april1}, combined with  a classical
one-dimensional Hardy inequality \cite[Theorem~5.1]{K}, ensure that
\[
\begin{split}
 \left(
 \int_0^{\infty}
  r^{\alpha p}
   \left|
     \frac{d^k}{dr^k}\left(\frac{f(r)}{r^{m-k}}\right)
   \right|^p
   dr
 \right)^{\frac{1}{p}}
& \le \frac{1}{(m-k-1)!}
 \left(
 \int_0^{\infty}
  r^{(\alpha -k-1)p}
\left(
  \int_0^r
s^k\left| \frac{d^m f}{ds^m} (s)\right|ds \right)^{p} dr
 \right)^{\frac{1}{p}}
\\
&
\le
C
 \left(
 \int_0^{\infty}
  r^{\alpha p}
   \left|
     \frac{d^m f}{dr^m} (r)
   \right|^p
   dr
 \right)^{\frac{1}{p}}
\end{split}
\]
for some
constant $C$. \qed

\bigskip
\par\noindent
The following result is a kind of combinatorial identity. In its
statement, we use the abridged notation
$$[k]=\{1,2,\ldots,k\}.$$

\begin{lemma}\label{lemma:2.3}
 Let $k\in \mathbb N$, and let
 $a_1,a_2,\ldots, a_k \in \mathbb R$. Given
 $I\subset [k]$, define $a_{\emptyset}=0$, and $a_{I}=\sum_{i\in I}a_{i}$ if $I\ne \emptyset$.
 Then
\begin{equation}
\label{eq:2.5}
 k!\prod^{k}_{i=1}a_i
 =
 \sum^{k}_{l=0}(-1)^{k-l}\sum_{I\subset [k] \atop \#
 I=l}(a_{I}+s)^{k} \quad \hbox{for  $s \in \mathbb R$,}
\end{equation}
 where $\# I$ denotes the cardinality of
$I$.
\end{lemma}

 \par\noindent {\bf Proof}. Define $\varphi_k : \mathbb R \to \mathbb R$ as 
\[
 \varphi_k(s) =\sum^{k}_{j=0}(-1)^{k-j}\sum_{I\subset
 [k]
 \atop \# I=j}(a_{I}+s)^{k} \quad \hbox{for $s \in \mathbb
 R$.}
\]
Equation \eqref{eq:2.5} can thus be rewritten as
\begin{equation}\label{april2}
 \varphi_k(s)= k!\prod^{k}_{i=1}a_i \quad \hbox{for $s \in \mathbb
 R$.}
 \end{equation}
We  establish equation \eqref{april2} by induction.
 The case when $k=1$ is trivial.
 Assume that \eqref{april2} holds for  $k-1$, for some $k\ge 2$.
 We begin by proving that $\varphi_k$ is constant.
 The following chain holds:
\begin{align*}
 \varphi_k(s)
&=\sum^{k}_{j=0}(-1)^{k-j}
\Biggl(\sum_{I\subset [k] \atop \# I=j, \,k\notin I}(a_{I}+s)^{k}+\sum_{I\subset [k] \atop \# I=j,\,k\in I}(a_{I}+s)^{k}\Biggr)\\
&=\sum^{k-1}_{j=0}(-1)^{k-j}\sum_{I\subset [k-1] \atop \#
I=j}(a_{I}+s)^{k}
+\sum^{k}_{j=1}(-1)^{k-j}\sum_{I\subset [k-1] \atop \# I=j-1}(a_{I}+(a_k+s))^{k}\\
&=-\sum^{k-1}_{j=0}(-1)^{(k-1)-j}\sum_{I\subset [k-1] \atop \#
I=j}(a_{I}+s)^{k}
+\sum^{k-1}_{j=0}(-1)^{k-(j-1)}\sum_{I\subset [k-1] \atop \# I=j}(a_{I}+(a_k+s))^{k}\\
&=-\sum^{k-1}_{j=0}(-1)^{(k-1)-j}\sum_{I\subset [k-1] \atop \#
I=j}(a_{I}+s)^{k} +\sum^{k-1}_{j=0}(-1)^{(k-1)-j}\sum_{I\subset
[k-1] \atop \# I=j}(a_{I}+(a_k+s))^{k}.
\end{align*}
Hence, by the induction assumption,
\begin{align*}
 \frac{d \varphi_{k}}{ds}(s)
&=-k\sum^{k-1}_{j=0}(-1)^{(k-1)-j}\sum_{I\subset [k-1] \atop \#
I=j}(a_{I}+s)^{k-1}
+k\sum^{k-1}_{j=0}(-1)^{(k-1)-j}\sum_{I\subset [k-1] \atop \# I=j}(a_{I}+(a_k+s))^{k-1}\\
&=-k\varphi_{k-1}(s)+k\varphi_{k-1}(a_k+s)=-k(k-1)!\prod^{k-1}_{i=1}a_i+k(k-1)!\prod^{k-1}_{i=1}a_i
=0 \quad \hbox{for $s \in \mathbb R$.}
\end{align*}
  Therefore $\varphi_k$ is constant. In particular,
\begin{equation}
\label{eq:7.1}
 \varphi_{k}(s)=\varphi_{k}(0) \quad \hbox{for $s \in \mathbb R$.}
\end{equation}
It remains to show that
\begin{equation}\label{april3}
\varphi_{k}(0) = k!\prod^{k}_{i=1}a_i.
\end{equation}
 To this purpose, observe that,  by the induction assumption,
\begin{equation}
\label{eq:7.2}
 k\int^{a_k}_{0}\varphi_{k-1}(s)ds
=k  (k-1)!\prod^{k-1}_{i=1}a_{i}\int^{a_k}_{0}ds
=k!\prod^{k}_{i=1}a_{i}.
\end{equation}
 On the other hand,  the very definition of $\varphi_{k-1}$ ensures
 that
\begin{align*}
\label{eq:7.1}
 k\int^{a_k}_{0}\varphi_{k-1}(s)ds
&=\int^{a_k}_{0}\sum^{k-1}_{j=0}(-1)^{k-1-j}\sum_{I\subset [k-1] \atop \# I=j}k(a_{I}+s)^{k-1}ds\\
&=\sum^{k-1}_{j=0}(-1)^{k-1-j}\sum_{I\subset [k-1] \atop \# I=j}\Bigl((a_{I}+a_k)^{k}-a_{I}^k\Bigr)\\
&=\sum^{k-1}_{j=0}(-1)^{k-(j+1)}\sum_{I\subset [k-1] \atop \#
I=j}(a_{I}+a_k)^{k}
-\sum^{k-1}_{j=0}(-1)^{(k-1)-j}\sum_{I\subset [k-1] \atop \# I=j}a_{I}^{k}\\
&=\sum^{k}_{j=1}(-1)^{k-j}\sum_{I\subset [k-1] \atop \#
I=j-1}(a_{I}+a_k)^{k}
+\sum^{k-1}_{j=0}(-1)^{k-j}\sum_{I\subset [k-1] \atop \# I=j}a_{I}^{k}\\
&=\sum^{k}_{j=0}(-1)^{k-j}\Biggl(\sum_{I\subset [k] \atop \#
I=j,\,k\in I}a_{I}^{k}+\sum_{I\subset [k] \atop \# I=j,k\notin
I}a_{I}^{k}\Biggl)
 =\sum^{k}_{j=0}(-1)^{k-j}\sum_{I\subset [k] \atop \# I=j}a_{I}^{k} =\varphi_{k}(0).
\end{align*}
 Hence,
\begin{equation}
\label{eq:7.3}
 k\int^{a_k}_{0}\varphi_{k-1}(s)ds=\varphi_{k}(0).
\end{equation}
 Equation \eqref{april3} follows from \eqref{eq:7.2} and
 \eqref{eq:7.3}.
 \qed

\bigskip
\noindent The next lemma is concerned with the special case when
$\alpha = \beta$, and hence $p=q$, in Theorem \ref{proposition:1.1}.

\begin{lemma}\label{lemma:2.2}
 Let $k,m\in \mathbb{N}$,   $m\ge 2$ and $1\le k\le m-1$.
Assume that $1\le p<\infty$ and $\alpha<k+\frac{p-1}{p}$. Then there exists a constant $C$ such
that
\begin{equation}\label{eq:2.3}
 \left(
\int_{\mathbb{R}^n_+}
  x_n^{\alpha p}
   \left|
    \nabla ^k \left(\frac{u(x)}{x_n^{m-k}}\right)
   \right|^p
   dx
 \right)^{\frac{1}{p}}
\le
C
 \left(
\int_{\mathbb{R}^n_+}
  x_n^{\alpha p}
   \left|
    \nabla ^m u
   \right|^p
   dx
 \right)^{\frac{1}{p}}
\end{equation}
for every $u\in C_0^{\infty}(\mathbb{R}^n_+)$.
\end{lemma}

\noindent {\bf Proof}.  The case when $n=1$ is the object of
Lemma~\ref{lemma:2.1}. We may thus assume that $n\ge 2$. Set
$\widetilde x = (x_1, \dots , x_{n-1})$, so that $x =
(\widetilde{x}, x_n)$ for $x \in \mathbb R^n$. Fubini's theorem and
Lemma~\ref{lemma:2.1} imply that
\begin{align}\label{eq:2.4}
  \left(
\int_{\mathbb{R}^n_+}
  x_n^{\alpha p}
   \left|
     \frac{\partial ^k}{\partial x_n^k}\left(\frac{u(x)}{x_n^{m-k}}\right)
   \right|^p
   dx
 \right)^{\frac{1}{p}}
& = \left( \int_{\mathbb{R}^{n-1}}\int_0^{\infty}
  x_n^{\alpha p}
   \left|
     \frac{\partial ^k}{\partial x_n^k}\left(\frac{u(\widetilde x, x_n)}{x_n^{m-k}}\right)
   \right|^p
   dx_n
   d\widetilde x
 \right)^{\frac{1}{p}}
\\ \nonumber
& \le C \left( \int_{\mathbb{R}^{n-1}}\int_0^{\infty}
  x_n^{\alpha p}
  \left|\frac{\partial ^m u}{\partial x_n^m} (\widetilde x, x_n)\right|^p dx_n d\widetilde x
 \right)^{\frac{1}{p}}
\\ \nonumber
& = C \left( \int_{\mathbb{R}_+^{n}}
  x_n^{\alpha p}
\left| \frac{\partial ^m  u}{\partial x_n^m}(x)\right|^p dx
 \right)^{\frac{1}{p}},
\end{align}
where $C$ denotes the constant appearing in inequality
\eqref{eq:2.1}. This establishes inequality \eqref{eq:2.3} with
$\nabla ^k$ replaced by the sole derivative $\frac{\partial
^k}{\partial x_n^k}$ on its left-hand side.
\\ Our next task is to extend this inequality to  arbitrary $k$-th order
derivatives. For ease of notation, let us set
$$\partial _\ell = \frac{\partial  }{\partial x_\ell}$$
for $\ell = 1, \dots , n$. Moreover, given $k$ indices $j_i\in
\{1,\cdots, n\}$, with $i=1, \dots ,k$, we define
$$  \partial^k
 =
 \prod_{i=1}^k \partial_{_{j_i}}.$$
  Lemma~\ref{lemma:2.3} implies
 that
\begin{equation}\label{eq:2.6}
k!\prod_{i=1}^k \partial_{_{j_i}}
=
 \sum^{k}_{l=0}(-1)^{k-l}\sum_{I\subset [k] \atop \# I=l}
 \Big(\sum_{i\in I}\partial_{j_i}+\partial_n\Big)^{k}
 .
\end{equation}
Now fix $I\subset [k]$. Then, there exist $h_l\in \{0,1,\cdots, n\}$
with $l=1, \dots , n$ and $\sum_{l=1}^nh_l=\# I$, such that
\[
 \sum_{i\in I}\partial_{j_i}+\partial_n
 =
  \sum_{l=1}^n h_l \partial_l  +\partial_n.
\]
Let $\Psi : \mathbb R^n \to \mathbb R^n$ be the bijective linear map
defined as
$$\Psi (y) = y + y_n \sum_{l=1}^n h_l  e_l \qquad \hbox{for $y \in \mathbb
R^n$,}
$$
where $\{e_1, \dots , e_n\}$ is the canonical  basis in
$\mathbb{R}^n$. Note that $\Psi (\mathbb R^n_+) = \mathbb R^n_+$,
and $\Psi (\partial \mathbb R^n_+) = \partial \mathbb R^n_+$.
Consider the change of variables
\[
 x=\Psi(y) \qquad \hbox{for $y \in \mathbb
R^n$.}
\]
Since its Jacobian determinant equals  $1+h_n$,
\begin{multline}\label{eq:2.7}
 \int_{\mathbb{R}^n_+}
  x_n^{\alpha p}
   \left|
\left(\sum_{l=1}^n h_l \partial_l  +\partial_n\right)^k
 \left(
  \frac{u(x)}{x^{m-k}_n}
 \right)
   \right|^p
   dx
\\
  = (1+h_n)\int_{\mathbb{R}^n_+}
  y_n^{\alpha p}
   \left|
\left( \sum_{l=1}^n h_l \partial_l  +\partial_n\right)^k \left(
  \frac{u(x)}{x^{m-k}_n}
 \right)_{ \lfloor {x=\Psi(y)}}
   \right|^p
   dy,
 \end{multline}
and
\begin{equation}\label{april4}
\int _{\mathbb{R}^n_+} x_n^{\alpha p} \left| \left( \sum_{l=1}^n h_l
\partial_l  +\partial_n\right)^{m}u(x)\right|^p dx = (1+ h_n)
\int _{\mathbb{R}^n_+} y_n^{\alpha p}\left| \left[ \left(
\sum_{l=1}^n h_l
\partial_l  +\partial_n\right)^{m}u(x)\right]
_{ \lfloor {x=\Psi(y)}}\right|^p\, dy
 \end{equation}
 for  $u\in C_0^{\infty}(\mathbb{R}^n_+)$.
Now define $v :  \mathbb R^n_+ \to  \mathbb R^n_+$ as
$v(y)=u(\Psi(y))$ for $y \in \mathbb R^n_+$. We have that $v \in
C^\infty _0(\mathbb R^n_+)$. By  the chain rule for derivatives and
the definition of $\Psi$,
\begin{align}\label{eq:2.8}
 \frac{\partial}{\partial y_n}
 \left(\frac{v(y)}{y_n^{m-k}}\right)
 & =
 \frac{\partial}{\partial y_n}
  \left(
 (1+h_n)^{m-k}
\left( \frac{u(x)}{x_n^{m-k}}
 \right)
 _{ \lfloor {x=\Psi(y)}}
 \right)
 \\ \nonumber
 &=
(1+h_n)^{m-k}  \left[ \left( \sum_{l=1}^n h_l
\partial_l  +\partial_n\right)
 \left(
  \frac{u(x)}{x^{m-k}_n}
 \right)
 \right]
_{ \lfloor {x=\Psi(y)}} \quad \hbox{for $y \in \mathbb R^n_+$.}
\end{align}
Iterating  equation \eqref{eq:2.8} yields
\begin{equation}\label{eq:2.9}
\begin{split}
 \frac{\partial^k}{\partial y_n^k}
 \left(\frac{v(y)}{y_n^{m-k}}\right)
 =
(1+h_n)^{m-k} \left[ \left( \sum_{l=1}^n h_l
\partial_l  +\partial_n\right)^k
 \left(
  \frac{u(x)}{x^{m-k}_n}
 \right)
\right]
 _{ \lfloor {x=\Psi(y)}}  \quad \hbox{for $y \in \mathbb R^n_+$.}
\end{split}
\end{equation}
Analogously,
\begin{equation}\label{eq:2.11}
 \frac{\partial^m v}{\partial y_n^m}(y)
 =
 \left[\left( \sum_{l=1}^n h_l \partial_l
 +\partial_n\right)^{m}u(x)\right]
 _{ \lfloor {x=\Psi(y)}}  \quad \hbox{for $y \in \mathbb R^n_+$.}
\end{equation}
Coupling  \eqref{eq:2.9} with \eqref{eq:2.4} tells us that
\begin{align}\label{eq:2.10}
& \left(\int_{\mathbb{R}^n_+}
  y_n^{\alpha p}
   \left|
  \left( \sum_{l=1}^n h_l \partial_l  +\partial_n\right)^k
 \left(
  \frac{u}{x^{m-k}_n}
 \right)
_{ \lfloor {x=\Psi(y)}}
   \right|^p
   dy\right)^{\frac 1p}
   \\ \nonumber
& \qquad \quad
 = (1+h_n)^{k-m} \left(
 \int_{\mathbb{R}^n_+}
y_n^{\alpha p}\left| \frac{\partial^k }{\partial
y_n^k}\left(\frac{v(y)}{y_n^{m-k}}\right)\right|^p
   dy
\right)^{\frac1p} \le C \left(
 \int_{\mathbb{R}^n_+}
y_n^{\alpha p}\left| \frac{\partial^m v}{\partial y_n^m} \right|^p
   dy
\right)^{\frac1p}.
\end{align}
Hence, via \eqref{eq:2.7}, \eqref{april4},  and \eqref{eq:2.11},
\begin{equation}\label{eq:2.12}
\left(
 \int_{\mathbb{R}^n_+}x_n^{\alpha p}
\left|
  \left(\sum_{i\in I}\partial_{j_i}+\partial_n\right)^{k}
 \left(
  \frac{u(x)}{x^{m-k}_n}
 \right)
 \right|^pdx
 \right)^{\frac1p}
\le C\left(
 \int_{\mathbb{R}^n_+}x_n^{\alpha p}
\left| \nabla ^m u \right|^pdx
 \right)^{\frac1p}
\end{equation}
for some constant $C>0$. Thanks to the arbitrariness of $I \subset
[k]$, inequality \eqref{eq:2.3} follows.
 \qed

\medskip
\par\noindent
Our last intermediate step consists in an estimate for the left-hand
side of inequality \eqref{eq:2.1}, involving $k$-th order
derivatives of $u(x)x_n^{k-m}$, in terms of the $k$-th and
$(k+1)$-th order derivatives of the same expression, but with
different weights.
Note that no
condition on the trace of $u$ on $\partial \mathbb R^n_+$ is
required here.

\begin{lemma}\label{lemma:2.4}  Let $k,m\in \mathbb{N}$,  $m\ge 2$ and $1\le k\le m-1$.
Assume  that
 $1\le p \le q <\infty$,
 $\frac1p-\frac1q= \frac{\beta-\alpha}{n}$,
 and $\alpha \le \beta \le \alpha+1$.
Then, there exists a constant $C$ such that
\begin{multline}\label{eq:2.13}
 \left(
  \int_{\mathbb{R}^n_+} x_n^{\beta q}
   \left|
    \nabla^k \left(\frac{u(x)}{x_n^{m-k}}\right)
   \right|^qdx
 \right)^{\frac1q}
 \\ \le
C
 \left(
  \int_{\mathbb{R}^n_+} x_n^{(\alpha+1)p}
   \left|
   \nabla ^{k+1}\left(\frac{ u(x)}{x_n^{m-k}}\right)
   \right|^p
   +
  x_n^{\alpha p}
   \left|
    \nabla^k \left( \frac{u(x)}{x_n^{m-k}}\right)
   \right|^p
   dx
 \right)^{\frac1p}
\end{multline}
for every $u \in C^\infty (\mathbb R_+^n)$.
\end{lemma}

%
\par\noindent {\bf Proof}.
By the same argument as in Lemma~\ref{lemma:2.2}, it suffices to
prove inequality \eqref{eq:2.13} with $\nabla ^k$ replaced with just
$\frac{\partial ^k}{\partial x_n^k}$ on the left-hand side. Let us
set, for simplicity of notation, $$\partial^k_n =
\frac{\partial ^k}{\partial x_n^k}.$$
 Assume first that  $1\le
p<n$.
  Let $\{B_j\}_{j\in \mathbb N}$  be a covering of
 $\mathbb{R}^n_+$
as in
 \cite[Lemma 4.1]{H1}, namely a family of balls $B_j$ with radius
 $r_j$, and centers in $\mathbb{R}^n_+$,
  such that:
\begin{itemize}
\item[\rm (i)]
$\mathbb{R}^n_+\subset \cup_{j=1}^{\infty}B_{j}$;
\item[\rm (ii)]
there exist positive constants $C'$ and $C''$ such that
 $
C' x_n \le r_j \le C'' x_n$ for every $j\in \mathbb{N}$ and $
 x\in B_j$;
\item[\rm (iii)]
there exists a positive constant $C$ such that $\#\{i\in
\mathbb{N}:B_j\cap B_i\neq \emptyset \}\le C$ for every $j\in
\mathbb{N}$.
\end{itemize}
 Let $\{\phi_j\}$ be a partition of unity subordinate to this
covering, namely a sequence of nonnegative functions $\phi _i \in
C^\infty _0 (B_j)$ such that $\sum _{j} \phi _j =1$.
The functions $\phi _j$ can be
chosen in such a way that $ |\nabla \phi_j|\le C/r_j $ for some
constant $C$, and for every $j \in \mathbb N$.
Set $p^*= \frac
{np}{n-p}$, the Sobolev conjugate of $p$. Since $\beta \leq \alpha +
1$, we have that $q \leq p^*$. Owing to the standard Sobolev
inequality for compactly supported functions, applied in each ball
$B_j$, the following chain holds:
\begin{align}\label{april10}
& \left( \int_{\mathbb{R}^n_+}
  x_n^{\beta q}
   \left|
    \partial^k_n \left(\frac{u(x)}{x_n^{m-k}}\right)
   \right|^{q}
   dx
 \right)^{\frac{1}{q}}
 \\ \nonumber & \qquad
\le \sum_{j\in \mathbb N}
 \left(
\int_{B_j}
  x_n^{\beta q}
   \left|
   \phi_j
    \partial^k_n \left(\frac{u(x)}{x_n^{m-k}}\right)
   \right|^{q}
   dx
 \right)^{\frac{1}{q}}
\le C_1\sum_{j\in \mathbb N} r_j^{\beta}
 \left(
\int_{B_j}
   \left|
   \phi_j
    \partial^k_n \left(\frac{u(x)}{x_n^{m-k}}\right)
   \right|^{q}
   dx
 \right)^{\frac{1}{q}}
 \\ \nonumber
  & \qquad
\le C_2\sum_{j\in \mathbb N} r_j^{\beta}
 \left(
\int_{B_j}
   \left|
   \phi_j
    \partial^k_n \left(\frac{u(x)}{x_n^{m-k}}\right)
   \right|^{p^*}
   dx
 \right)^{\frac{1}{p^*}}
 \left(
\int_{B_j}
1dx
 \right)^{\frac{1}{q}-\frac{1}{p^*}}
 \\ \nonumber & \qquad \le C_3\sum_{j\in \mathbb N} r_j^{1 +\alpha} \left(\int_{B_j}
   \left|
   \nabla
   \left(
   \phi_j
    \partial^k_n \left(\frac{u(x)}{x_n^{m-k}}\right)
    \right)
   \right|^p
   dx
 \right)^{\frac 1p}
 \\ \nonumber
 & \qquad
\le C_4 \left( \int_{\mathbb{R}^n_+} x_n^{\alpha p}
   \left|
   \left(
    \partial^k_n \left(\frac{u(x)}{x_n^{m-k}}\right)
    \right)
   \right|^p
   dx
\right)^{\frac 1p} + C_4 \left( \int_{\mathbb{R}^n_+}
x_n^{(\alpha+1)p}
   \left|
   \nabla
   \left(
    \partial^k_n \left(\frac{u(x)}{x_n^{m-k}}\right)
    \right)
   \right|^p
   dx
 \right)^{\frac{1}{p}},
\end{align}
for suitable constants $C_i$, $i=1, \dots 4$, for every $u \in
C^\infty(\mathbb R^n_+)$. This establishes inequality
\eqref{eq:2.13} in the case when $1 \leq p <n$.  Note that, in fact,
this argument  proves \eqref{eq:2.13}
for every $(k+1)$-times weakly differentiable function $u$ in
$\mathbb{R}^n_+$ making the right-hand side finite.
\par\noindent
Assume next that $n\le p <\infty$.  Let $p_1, q_1,\alpha_1,\beta_1,
r$ be such that
\[
\begin{split}
\max\left\{1, n\left(1-\frac{p}{q}\right)\right\}\le p_1<n, \quad
 r=\frac{p}{p_1}
 ,\quad q_1=\frac{q}{r},\quad \alpha_1=r \alpha, \quad
 \beta_1= r\beta .
\end{split}
\]
Hence,
\begin{equation}\label{eq:2.14}
\begin{cases}
rp_1=p,\quad rq_1=q,\quad \alpha_1p_1=\alpha p, \quad
\beta_1q_1=\beta q,
\quad
p_1(\alpha_1+r)=p(\alpha+1),
\\
 \frac{\beta_1-\alpha_1}{n}
 =
 \frac{r}{n}(\beta-\alpha)
 =
 r\big(\frac1p-\frac1q\big)
 =
 \frac{1}{p_1}-\frac{1}{q_1},
 \\
 0\le \beta_1-\alpha_1=n\big(\frac{1}{p_1}-\frac{1}{q_1}\big)
 =\frac{n}{p_1}\big(1-\frac{p}{q}\big)\le 1.
 \end{cases}
\end{equation}
Given any function $u\in C^\infty (\mathbb{R}^n_+)$,   define $w :
\mathbb R^n_+ \to \mathbb R$ as
\[
 w(\widetilde x,x_n)
 =
 x_n^{m-k}
 \int_0^{x_n}
 \int_0^{s_{k-1}}
\cdots
\int_0^{s_1}
  \left|
   \partial_n^k
    \left(
     \frac{u(\widetilde x,s)}{s^{m-k}}
    \right)
  \right|^rdsds_1\cdots ds_{k-1}
\]
for $(\widetilde x, x_n) \in \mathbb R^n_+$. Note that the function
$w$ is $k+1$-times weakly differentiable. An application of
inequality
 \eqref{april10} to $w$, with $\alpha , \beta, p, q$ replaced by  $\alpha _1 , \beta _1, p_1,
 q_1$, yields
\begin{multline}\label{eq:2.15}
 \left(
  \int_{\mathbb{R}^n_+} x_n^{\beta_1 q_1}
   \left|
    \partial_n^k\left(\frac{w(x)}{x_n^{m-k}}\right)
   \right|^{q_1}dx
 \right)^{\frac1{q_1}}
 \\
  \le
C
 \left(
  \int_{\mathbb{R}^n_+} x_n^{(\alpha_1+1)p_1}
   \left|
   \nabla  \partial_n^k\left(\frac{w(x)}{x_n^{m-k}}\right)
   \right|^{p_1}
   dx
 \right)^{\frac1{p_1}}
   +
C \left( \int_{\mathbb{R}_+^n}  x_n^{\alpha_1 p_1}
   \left|
    \partial_n^k \left( \frac{w(x)}{x_n^{m-k}}\right)
   \right|^{p_1}
   dx
 \right)^{\frac1{p_1}}.
\end{multline}
By the very definition of $w$,
\begin{equation}\label{eq:2.16}
 \left(
\int_{\mathbb{R}^n_+}
  x_n^{\beta_1 q_1}
   \left|
    \partial_n^k \left(\frac{w(x)}{x_n^{m-k}}\right)
   \right|^{q_1}
   dx
 \right)^{\frac{1}{q_1}}
=
 \left(
\int_{\mathbb{R}^n_+}
  x_n^{\beta_1 q_1}
   \left|
    \partial_n^k \left(\frac{u(x)}{x_n^{m-k}}\right)
   \right|^{rq_1}
   dx
 \right)^{\frac{1}{q_1}}.
\end{equation}
Moreover,
\begin{equation}\label{eq:2.17}
\begin{split}
\left|
 \nabla \partial_n^k\left(\frac{w(x)}{x_n^{m-k}}\right)
\right| = \left|
 \nabla \left|\partial_n^k\left(\frac{u(x)}{x_n^{m-k}}\right)\right|^r
\right| =r \left| \nabla \partial_n^k
\left(\frac{u(x)}{x_n^{m-k}}\right) \right| \left|\partial_n^k
\left(\frac{u(x)}{x_n^{m-k}}\right)\right|^{r-1}
\end{split}
\end{equation}
for $x \in \mathbb{R}^n_+$. From equation  \eqref{eq:2.17}  and
Young's inequality one can infer that
\begin{multline}\label{eq:2.18}
\left(
  \int_{\mathbb{R}^n_+} x_n^{(\alpha_1+1)p_1}
   \left|
   \nabla   \partial_n^k \left(\frac{w(x)}{x_n^{m-k}}\right)
   \right|^{p_1}
   dx
\right)^{\frac{1}{p_1}}
\\
\le C \left( \int_{\mathbb{R}^n_+} x_n^{\alpha_1 p_1}
 \left|
  \partial_n^k\left(\frac{u(x)}{x_n^{m-k}}\right)
 \right|^{rp_1}
dx \right)^{\frac{1}{p_1}} +C \left( \int_{\mathbb{R}^n_+}
x_n^{(\alpha_1+r)p_1}
 \left|
  \nabla \partial_n^k \left(\frac{u(x)}{x_n^{m-k}}\right)
 \right|^{rp_1}
   dx
 \right)^{\frac{1}{p_1}}
\end{multline}
for some constant $C$. Combining  \eqref{eq:2.15}, \eqref{eq:2.16},
and \eqref{eq:2.18}, we have, owing to the first line in
\eqref{eq:2.14},
\[
\begin{split}
& \left( \int_{\mathbb{R}^n_+}
  x_n^{\beta q}
   \left|
    \partial_n^k \left(\frac{u(x)}{x_n^{m-k}}\right)
   \right|^{q}
   dx
 \right)^{\frac{1}{q}}
\\
& \quad \le C \left( \int_{\mathbb{R}^n_+} x_n^{\alpha p}
 \left|
  \partial_n^k\left(\frac{u(x)}{x_n^{m-k}}\right)
 \right|^{p}
dx \right)^{\frac{1}{p}} +C \left( \int_{\mathbb{R}^n_+}
x_n^{(\alpha+1)p}
 \left|
  \nabla \partial_n^k\left(\frac{u(x)}{x_n^{m-k}}\right)
 \right|^{p}
   dx
 \right)^{\frac{1}{p}},
\end{split}
\]
namely \eqref{eq:2.13}. The proof is complete. \qed

\medskip
\par\noindent
We are now in a position to accomplish the proof of Theorem
\ref{proposition:1.1}.

\medskip
\par\noindent
{\bf Proof of Theorem~\ref{proposition:1.1}}. Let $p,q,\alpha,\beta$
be as in the statement. For $i=0,\cdots, m-k$, define $p_i,\alpha_i$
as:
\[
 \alpha_i
 =\frac{(m-k-i)\beta+i\alpha}{m-k},
\qquad
 \frac{1}{p_i}
 =
 \frac{1}{m-k}
  \left(
   \frac{i}{p}+\frac{m-k-i}{q}
  \right).
\]
One can verify that
\begin{equation}\label{eq:2.19}
\begin{split}
1\le p_{i+1}\le p_i <\infty,\quad
\frac{1}{p_{i+1}}-\frac{1}{p_i}=\frac{\alpha_i-\alpha_{i+1}}{n},\quad
\alpha_i-\alpha_{i+1}=\frac{\beta-\alpha}{m-k}\le 1,
\end{split}
\end{equation}
for every $i=0,\cdots, m-k$, and
$$ \alpha _0 = \beta, \quad \alpha _{m-k} = \alpha, \quad p_0 =q, \quad
p_{m-k}=p\,.$$ An iterated application of Lemma~\ref{lemma:2.4},
with $\alpha=\alpha_{i+1},\  \beta=\alpha_i,\ q=p_i,\  p=p_{i+1}$,
yields:
\begin{align}\label{eq:2.20}
 & \left(
  \int_{\mathbb{R}^n_+}
  x_n^{\beta q}
   \left|
    \nabla^k \left(\frac{u(x)}{x_n^{m-k}}\right)
   \right|^qdx
 \right)^{\frac1q}
\\ \nonumber
&
\le
C
 \left(
  \int_{\mathbb{R}^n_+}
  x_n^{ (\alpha_1+1)p_1}
   \left|
  \nabla^{k+1} \left(\frac{u(x)}{x_n^{m-k}}\right)
   \right|^{p_1}dx
 \right)^{\frac{1}{p_1}}
+
 \left(
  \int_{\mathbb{R}^n_+}
  x_n^{ \alpha_1 p_1}
   \left|
    \nabla^k \left(\frac{u(x)}{x_n^{m-k}}\right)
   \right|^{p_1}dx
 \right)^{\frac{1}{p_1}}
\\ \nonumber
&
\le
C
\sum_{j=0}^2
 \left(
  \int_{\mathbb{R}^n_+}
  x_n^{ (\alpha_2+j)p_2}
   \left|
     \nabla^{k+j} \left(\frac{u(x)}{x_n^{m-k}}\right)
   \right|^{p_2}dx
 \right)^{\frac{1}{p_2}} 
\\ \nonumber & \qquad \vdots
\\ \nonumber
&
\le
C
\sum_{j=0}^{m-k}
 \left(
  \int_{\mathbb{R}^n_+}
  x_n^{ (\alpha+j)p}
   \left|
    \nabla ^{k+j} \left(\frac{u(x)}{x_n^{m-k}}\right)
   \right|^{p}dx
 \right)^{\frac{1}{p}},
\end{align}
for every $u \in C^\infty _0(\mathbb R^n_+)$.
If we show that
\begin{equation}\label{eq:2.21}
\left|
  \nabla^{k+j} \left(\frac{u(x)}{x_n^{m-k}}\right)
\right|
 \le
  \frac{C}{x_n^j}\sum_{i=k}^{k+j}\left|
  \nabla^i\left(\frac{u(x)}{x_n^{m-i}}\right)\right| \quad \hbox{for
  $x \in \mathbb R^n_+$,}
\end{equation}
for some constant $C>0$, and for any such function $u$, then the
conclusion will follow via Lemma~\ref{lemma:2.2} and inequality
\eqref{eq:2.20}.
\\ In order to prove \eqref{eq:2.21}, let us denote by
$\widetilde\partial ^i$
any partial derivative of order $i$ with respect to the
variables  $x_1,\dots, x_{n-1}$. Then
any derivative in $\nabla ^{k+j}$
can be written as $\widetilde\partial^{k+j-l}\partial_n^{l}$ for
some $l=0,1,\cdots, k+j$. Inequality \eqref{eq:2.21} trivially holds
if $l =0$. Assume now that $l\ge 1$. Let us preliminarily observe
that
\begin{equation}\label{eq:2.22}
\partial_n^{l}\left(\frac{u(x)}{x_n^{m-k}}\right)
= \frac{1}{x_n} \left[
 \partial_n^{l}
 \left(\frac{u(x)}{x_n^{m-k-1}}\right)
 -l
 \partial_n^{l-1}\left(\frac{u(x)}{x_n^{m-k}}\right)
\right] \quad \hbox{for
  $x \in \mathbb R^n_+$.}
\end{equation}
Equation \eqref{eq:2.22} can be verified by induction on $l$. The
case when $l=1$ is easy, since
\[
 \partial_n
 \left(\frac{u(x)}{x_n^{m-k}}\right)
= \frac{1}{x_n}\left[
\partial_n
 \left(\frac{u(x)}{x_n^{m-k-1}}\right)
- \frac{u(x)}{x_n^{m-k}}
 \right] \quad \hbox{for
  $x \in \mathbb R^n_+$.}
\]
 Assume next that \eqref{eq:2.22} holds with $l$ replaced by $l-1$, for some  $l\ge 2$.
Computations show  that
\[
 \partial_n\left[x_n\partial_n^{l-1}\left(\frac{u(x)}{x_n^{m-k}}\right)\right]
 =\partial_n^{l-1}\left(\frac{u(x)}{x_n^{m-k}}\right)+x_n\partial_n^{l}\left(\frac{u(x)}{x_n^{m-k}}\right) \quad \hbox{for
  $x \in \mathbb R^n_+$.}
\]
On the other hand, the induction assumption ensures that
\[
  \partial_n\left[x_n\partial_n^{l-1}\left(\frac{u(x)}{x_n^{m-k}}\right)\right]
=
 \partial_n^{l}
 \left(\frac{u(x)}{x_n^{m-k-1}}\right)
 -(l-1)
 \partial_n^{l-1}\left(\frac{u(x)}{x_n^{m-k}}\right) \quad \hbox{for
  $x \in \mathbb R^n_+$.}
\]
Equation \eqref{eq:2.22} is thus   established   for every $l\ge 1$.
\\
Iterating equation \eqref{eq:2.22} $j$ times tells us that
$$
 \partial_n^{l}\left(\frac{u(x)}{x_n^{m-k}}\right)
 =
\frac{1}{x_n^j}\sum_{i=\max\{0, l-j\}}^{l}c_i
\partial_n^i\left(\frac{u(x)}{x_n^{m-k+l-i-j}}\right) \quad \hbox{for
  $x \in \mathbb R^n_+$,}
%
%
%
%
$$
for some constants $c_i \in \mathbb{R}$. Therefore,
\begin{align*}
\left|
 \widetilde{\partial}^{j+k-l}\partial_n^l\left(\frac{u(x)}{x_n^{m-k}}\right)
\right| &
 \leq C
\frac{1}{x_n^j}\sum_{i=\max\{0,
l-j\}}^{l}\left|\widetilde{\partial}^{j+k-l}\partial_n^i\left(\frac{u(x)}{x_n^{m-k+l-i-j}}\right)\right|
\\ &\leq
\frac{C}{x_n^j}\sum_{h=k}^{k+j}\left|\nabla
^h\left(\frac{u(x)}{x_n^{m-h}}\right)\right| \quad \hbox{for
  $x \in \mathbb R^n_+$,}
\end{align*}
%
%
for some constant $C>0$. This establishes  inequality
\eqref{eq:2.21}. \qed

\section{Proof of the main result}\label{bounded}

Given a bounded smooth open set $\Omega$ with smooth boundary in
$\mathbb R^n$, with
$n \geq 1$,
we
make use of an orthogonal coordinate system which, in a sense,
rectifies $\partial \Omega$ in a suitable  neighborhood inside
$\Omega$. By the latter expression, we mean a subset $\Omega
_\varepsilon$ of $\Omega$ of the form
\begin{equation}\label{omegaeps}
 \Omega_{\varepsilon}=\{x\in \Omega: d(x)<\varepsilon\}
\end{equation}
for some $\varepsilon>0$. 
Let $\varepsilon_0$ be small enough for $d$ to agree, in
$\Omega_{\varepsilon_0}$,  with the distance function from $\partial
\Omega$.
 It is well known
that $\varepsilon_0$ can be chosen so small that, for every $x\in
\Omega_{\varepsilon_0}$, there exists a unique $y_x\in
\partial \Omega$ fulfilling
\begin{equation}\label{xy}
x=y_x+
d(x)\nu (y_x),
\end{equation}
where $\nu$ denotes  the  inward unit normal vector on $\partial
\Omega$.
\\
 Since $\partial \Omega$ is smooth, for every $x_0\in \partial
 \Omega$ there exist an open
neighborhood $\mathcal{U} (x_0)$ of $x_0$ on $ \partial \Omega$, a
radius ${r_0}>0$, and a smooth diffeomorphism
$$
\widetilde\Phi: B^{n-1}_{r_0}(0)\to \mathcal{U}(x_0).
$$
Next, define
 $\Phi: B^{n-1}_{r_0}(0) \times(0 ,\varepsilon_0)\to
\mathbb{R}^n$ as
\[
 \Phi(y)=\widetilde{\Phi}(\widetilde y) + y_n\nu (\widetilde \Phi(\widetilde y)) \quad \hbox{for $y \in B^{n-1}_{r_0}(0) \times(0
 ,\varepsilon_0)$,}
\]
where $y = (\widetilde y, y_n)$, and $\widetilde y =(y_1,\cdots,
y_{n-1})$.
By \eqref{xy},
\begin{equation}\label{ostrov10}
y_n = d(\Phi (y)) \quad \hbox{for $y \in B^{n-1}_{r_0}(0) \times(0
 ,\varepsilon_0)$.}
 \end{equation}
On setting
\begin{equation}\label{N}
\mathcal{N}(x_0) = \{
 x\in \Omega_{\varepsilon_0}:y_x\in \mathcal{U}(x_0)
\},
\end{equation}
one can prove that the map $\Phi : {B^{n-1}_{r_0} (0)
\times(0,\varepsilon_0)} \to \mathcal{N}(x_0)$ is a smooth
diffeomorphism. In particular,
\[
 \mathcal{N}(x_0)=\Phi(B^{n-1}_{r_0}(0) \times (0,\varepsilon_0)).
\]
As a consequence, there exists a positive constant $C$ such that
\begin{equation}\label{eq:4.1}
 \frac{1}{C}\int_{B^{n-1}_{r_0}(0)}\int_0^{\varepsilon_0}|f(\Phi(y))|dy_nd\widetilde y
\le \int_{\mathcal{N}(x_0)}|f(x)|dx \le
C\int_{B^{n-1}_{r_0}(0)}\int_0^{\varepsilon_0}|f(\Phi(y))|dy_nd\widetilde
y
\end{equation}
for every function $f\in L^1(\mathcal{N}(x_0))$.

\medskip
\par\noindent
 In preparation for the proof of Theorem \ref{proposition:1.3},
we establish  the following local version.

\begin{lemma}\label{lemma:4.1}
Let $\Omega$, $p, q, r, m, k$ be as in the statement of Theorem
\ref{proposition:1.3}. Given any point $x_0\in
\partial \Omega$, let $\mathcal{N}(x_0)$ be defined as  in \eqref{N}.
 Then there exists a constant $C$ such that
\begin{equation}\label{april30}
\left(
 \int_{\mathcal{N}(x_0)}
  d(x)^{r}
   \left|
    \nabla ^k
     \left(
      \frac{u(x)}{d(x)^{m-k}}
     \right)
   \right|^q
  dx
\right)^{\frac1q}
\le
C
\sum_{l=1}^{m}
\left(
 \int_{\mathcal{N}(x_0)}
  d(x)^{p-1}
    |
    \nabla ^l
    u
   |^p
  dx
\right)^{\frac1p}
\end{equation}
for every $u\in C_0^{\infty}(\mathcal{N}(x_0))$.
\end{lemma}

\smallskip
\par\noindent
{\bf Proof}. Fix any function $u\in C_0^{\infty}(\mathcal{N}(x_0))$.
Owing to inequality \eqref{eq:4.1},
\begin{equation}\label{eq:4.2}
 \left(
 \int_{\mathcal{N}(x_0)}
  d(x)^{r}
   \left|
    \nabla^k
     \left(
      \frac{u(x)}{d(x)^{m-k}}
     \right)
   \right|^q
  dx
\right)^{\frac1q}
\le
C
\left(
 \int_{B^{n-1}_{r_0}(0)}\int_0^{\varepsilon_0}
  y_n^{r}
   \left|
   \nabla^k
     \left(
      \frac
      {u(x)}
      {d(x)^{m-k}}
     \right)_{
\lfloor {x=\Phi(y)}}
   \right|^q
  d\widetilde y dy_n
\right)^{\frac1q}
\end{equation}
for some constant $C>0$. Define
$$
 v(y) =u(\Phi(y)) \quad \hbox{and} \quad
 \delta (y) =d(\Phi(y)) \qquad \hbox{for $y \in B^{n-1}_{r_0}(0) \times
 (0,\varepsilon_0)$.}
$$
Observe that, by \eqref{ostrov10},
 $\delta(y)=y_n$ for $y \in
B^{n-1}_{r_0}(0) \times (0,\varepsilon_0)$. By the  chain
rule for derivatives,
\begin{equation}\label{eq:4.3}
    \left|
   \nabla^k
     \left(
      \frac{u(x)}{d(x)^{m-k}}
     \right)_{
\lfloor {x=\Phi(y)}}
   \right|
\le C \sum_{l=1}^k \left|  \nabla^l
     \left(
      \frac{v(y)}{\delta(y)^{m-k}}
     \right)
 \right|
  \qquad \hbox{for $y \in B^{n-1}_{r_0}(0) \times
 (0,\varepsilon_0)$,}
\end{equation}
for some constant $C$.
 Inequality
\eqref{eq:4.3} implies that
\begin{multline}\label{eq:4.4}
\left(
 \int_{B^{n-1}_{r_0}(0)}\int_0^{\varepsilon_0}
  y_n^{r}
   \left|
   \nabla^k
     \left(
      \frac{u(x)}{d(x)^{m-k}}
     \right)_{
\lfloor_{x=\Phi(y)}}
   \right|^q
  d\widetilde y dy_n
\right)^{\frac1q}
\\
 \le C \sum_{l=1}^k \left(
 \int_{B^{n-1}_{r_0}(0)}\int_0^{\varepsilon_0}
  y_n^{r}
   \left|
  \nabla^l
     \left(
      \frac{v(y)}{\delta(y)^{m-k}}
     \right)
   \right|^q
  d\widetilde y dy_n
\right)^{\frac1q}.
\end{multline}
It follows from Theorem~\ref{proposition:1.1},
applied with $m$ replaced with $m-k+l$, and with $k=l$, $\alpha =
\frac{p-1}p$, and $\beta=\frac{n}{p}-\frac{n}{q}+\frac{p-1}{p}$,
that
\begin{align}
\label{eq:4.5} & \left(
 \int_{B^{n-1}_{r_0}(0)}\int_0^{\varepsilon_0}
  y_n^{r}
   \left|
   \nabla^l
     \left(
      \frac{v(y)}{\delta(y)^{m-k}}
     \right)
   \right|^q
  d\widetilde y dy_n
\right)^{\frac1q}
\\ \nonumber
& \qquad \qquad \le \left(\varepsilon_0^{r-\beta q}
 \int_{B^{n-1}_{r_0}(0)}\int_0^{\varepsilon_0}
  y_n^{\beta q}
   \left|
   \nabla^l
     \left(
      \frac{v(y)}{\delta(y)^{m-k}}
     \right)
   \right|^q
  d\widetilde y dy_n
\right)^{\frac1q}
\\
 \nonumber &
\qquad
\qquad
\le C \left(
 \int_{B^{n-1}_{r_0}(0)}\int_0^{\varepsilon_0}
  y_n^{p-1}
   |
\nabla ^{m-k+l}
      {v(y)}
   |^p
  d\widetilde y dy_n
\right)^{\frac1p}
\end{align}
for some positive constant $C$.
Observe that the first inequality
holds owing to condition \eqref{r}.
The chain rule again ensures
that
\[
    |\nabla^{m-k+l}
      {v(y)}|
 \le
 C\sum_{j=1}^{m-k+l}
  |\nabla^j u(x)|_{\lfloor {x=\Phi(y)}} \qquad \hbox{for $y \in B^{n-1}_{r_0}(0) \times
 (0,\varepsilon_0)$.}
\]
 Hence, by \eqref{eq:4.1} and \eqref{ostrov10},
\begin{equation}\label{eq:4.6}
\sum_{l=1}^{k}
\left(
 \int_{B^{n-1}_{r_0}(0)}\int_0^{\varepsilon_0}
  y_n^{p-1}
   |
\nabla^{m-k+l}
      {v(y)}
   |^p
  d\widetilde y dy_n
\right)^{\frac 1p} \le C \sum_{j=1}^{m} \left(
 \int_{\mathcal{N}(x_0)}
  d(x)^{p-1}
   \left|
\nabla^{j}
      {u(x)}
   \right|^p
  dx
\right)^{\frac1p}
\end{equation}
for some constant $C$. The conclusion follows from inequalities
\eqref{eq:4.2}, \eqref{eq:4.4}, \eqref{eq:4.5}, and \eqref{eq:4.6}.
\qed

\bigskip

%
%

\par\noindent {\bf Proof of Theorem~\ref{proposition:1.3}}. Without loss of generality,
we can assume that $u\in C_0^{\infty}(\Omega)$.  Since $\partial
\Omega$ is compact and $\{\mathcal{U}(x):x\in \partial\Omega\}$ is
an open covering of $\partial \Omega$, there exist $N\in \mathbb{N}$
and $\{x_l\}_{l=1}^{N}\subset \partial \Omega$ such that $
\displaystyle \Omega _{\varepsilon _0} =
 \cup_{l=1}^{N}\mathcal{N}(x_l)$, where $\Omega _{\varepsilon
 _0}$ and $\mathcal{N}(x_l)$ are defined as in \eqref{omegaeps} and
 \eqref{N}, and $\varepsilon _0$ is chosen in such a way that
 \eqref{xy} holds.
Let $\{\phi_l\}_{l=0}^{N}$ be a partition of unity of functions
$\phi_l \in C_0^{\infty}(\mathbb{R}^n)$ such that
\begin{itemize}
\item[(i)] $0\le \phi_l\le 1$ for   $l=0,\cdots, N$, and $\sum_{l=0}^{N}\phi(x)=1$ for   $x\in
\Omega$;
\item[(ii)] $\text{supp}\, \phi_l  \cap \Omega \subset  \mathcal{N}(x_l)$ for
$l=1,\cdots,N$;
\item[(iii)] $\text{supp} \,\phi_0 \subset \Omega$.
\end{itemize}
Set $u_l =\phi_l u$. Then
\begin{multline}\label{eq:4.7}
\left(
 \int_{\Omega}d(x) ^{r}
  \left|
   \nabla^k\left(\frac{u (x)}{d(x) ^{m-k}}\right)
  \right|^qdx
\right)^{\frac1q}
\\
  \le C \sum_{l=1}^N \left(
 \int_{\mathcal{N}(x_l)}d(x)^{r}
  \left|
   \nabla^k\left(\frac{u_l(x)}{d(x)^{m-k}}\right)
  \right|^qdx
\right)^{\frac{1}{q}}
+
\left(
 \int_{\text{supp}\,\phi_0}d(x)^{r}
  \left|
   \nabla^k\left(\frac{u_0(x)}{d(x)^{m-k}}\right)
  \right|^qdx
\right)^{\frac{1}{q}}.
\end{multline}
Observe that there exists a positive constant $C$ such that $C\le
d(x)\le 1/C>0$, and $|\nabla ^h d(x) (x)|\leq C$ for $h=1, \dots ,
m$ and for $x\in \text{supp}\,\phi_0$. Thus, owing to assumption
\eqref{pq}, the standard Sobolev inequality ensures that
%
%
\begin{equation}\label{eq:4.8}
\left(
 \int_{\text{supp}\,\phi_0}d(x)^{r}
  \left|
   \nabla^k\left(\frac{ u_0(x)}{d(x)^{m-k}}\right)
  \right|^qdx
\right)^{\frac{1}{q}} \le C \sum_{j=0}^{ m} \left(
 \int_{\text{supp}\,\phi_0}
 d(x)^{p-1}
  |
   \nabla^ju_0
  |^pdx
\right)^{\frac{1}{p}}
\end{equation}
for some
constant $C$. On the other hand, Lemma~\ref{lemma:4.1}
tells us that
\begin{equation}\label{eq:4.9}
 \left(
 \int_{\mathcal{N}(x_l)}d(x)^{r}
  \left|
   \nabla^k\left(\frac{ u_l (x)}{d(x)^{m-k}}\right)
  \right|^qdx
\right)^{\frac{1}{q}}
\le
C\sum_{j=1}^m
\left(
 \int_{\mathcal{N}(x_l)}d(x)^{p-1}
   |
   \nabla ^j u_l
   |^pdx
\right)^{\frac{1}{p}}
\end{equation}
for $l=1 , \dots , N$. Inequalities \eqref{eq:4.7}, \eqref{eq:4.8},
and \eqref{eq:4.9} yield
\begin{align}
&
\left(
 \int_{\Omega}d(x)^{r}
  \left|
   \nabla^k\left(\frac{u (x)}{d(x)^{m-k}}\right)
  \right|^qdx
\right)^{\frac1q}
\\ \nonumber
&  \le C\sum_{l=1}^N\sum_{j=1}^m \left(
 \int_{\mathcal{N}(x_l)}d(x)^{p-1}
  \left|
   \nabla^j u_l
  \right|^pdx
\right)^{\frac{1}{p}} + C
\sum_{j=0}^{
m
} \left(
 \int_{\text{supp}\,\phi_0}
 d(x)^{p-1}
   |
   \nabla ^j
   u_0
   |^pdx
\right)^{\frac{1}{p}}
\\ \nonumber
&  \le C' \sum_{j=0}^m \left(
 \int_{\Omega}d(x)^{p-1}
  |
   \nabla^j u
  |^pdx
\right)^{\frac1p}.
\end{align}
It  remains to show that there exists a constant $C$ such that
\begin{equation}\label{april40}
\left(
 \int_{\Omega}d(x)^{r}
  \left|
   \nabla ^h\left(\frac{u(x)}{d(x)^{m-k}}\right)
  \right|^qdx
\right)^{\frac1q}
 \le C \sum_{j=0}^m \left(
 \int_{\Omega}d(x)^{p-1}
  |
   \nabla^j u
  |^pdx
\right)^{\frac1p}
\end{equation}
for $h=0, \dots , k-1$ as well. To this purpose,  we apply a
Hardy-Sobolev inequality  from \cite[Theorem~3]{H1}, which tells us
what follows. Assume that $1\le  p\le  q <\infty$ and
$\widetilde \alpha , \widetilde
\beta \in \mathbb R$
fulfil the conditions:
\begin{equation}\label{condition:horiuchi}
-1/{ q} < \widetilde \beta \le \widetilde \alpha, \qquad
m-k+1
-\widetilde \alpha +\widetilde
\beta<\frac{n}{ p}\le
m -k +1
-\widetilde \alpha +\widetilde \beta+\frac{n}{q}.
\end{equation}
Then there exists a constant $C>0$ such that
\begin{equation}\label{eq:horiuchi}
\|u\|_{W^{k-1,  q}(\Omega, d^{
\widetilde \beta q})} \le C
\|u\|_{W^{m,  p}
(\Omega,d^{
\widetilde \alpha p})}
\end{equation}
for every $u\in  W^{m,  p}(\Omega,d^{   \widetilde \alpha p})$.
Choose
\[
\widetilde \beta=\frac{r}{q},
\]
and $\widetilde \alpha$ such that
\begin{equation}\label{ostrov}
 \max\left\{\frac{p-1}{p}+m-k,\frac{r}{q}-\frac{n}{p}+m-k+1,\frac{r}{q}\right\}<\widetilde \alpha \le m-k+1+\frac{r}{q}-\frac{n}{p}+\frac{n}{q},\quad
 \widetilde \alpha \notin \mathbb{N}-\{1/p\}.
\end{equation}
 Conditions
\eqref{condition:horiuchi} are fulfilled with such a  choice of
$\widetilde \alpha$ and $\widetilde \beta$, which is possible thanks
to assumptions \eqref{pq} and \eqref{r}. Consequently,
\[
 \left\|
  \frac{u}{d^{m-k}}
 \right\|_{W^{k-1,q}(\Omega,d^r)}
 \le
 C
 \left\|
  \frac{u}{d^{m-k}}
 \right\|_{W^{m,p}(\Omega,d^{\widetilde \alpha p})}
 \le
 C
\sum_{j=0}^m\sum_{i=0}^j
 \left\|
  \frac{\nabla ^iu}{d^{m-k+j-i}}
 \right\|_{L^{p}(\Omega,d^{\widetilde \alpha p})}
\]
for every $u\in C_0^{\infty}(\Omega)$. Since $\widetilde \alpha
\notin \mathbb{N}-\frac{1}{p}$, we have that $p(\widetilde
\alpha-(m-k+j-i))\neq -1$ if $0\le j\le m,\ 0\le i\le j$. Hence, by
a standard  Hardy-Sobolev  embedding
\cite[Equation (8.37)]{K},
\[
  \sum_{i=0}^j
 \left\|
  \frac{\nabla ^iu}{d^{m-k+j-i}}
 \right\|_{L^{p}(\Omega,d^{\widetilde \alpha p})}
 \le
 C
 \sum _{l=0}^j
 \left\|
  \frac{\nabla^l u}{d^{m-k}}
 \right\|_{L^{p}(\Omega,d^{\widetilde \alpha p})}
\]
for some constant $C>0$, and for every function $u\in
C_0^{\infty}(\Omega)$. Inasmuch as the function $d$ is bounded in
$\Omega$ and, by \eqref{ostrov}, $\frac{p-1}{p}+m-k <\widetilde
\alpha$, one hence deduces that
\[
 \left\|
  \frac{\nabla^ju}{d^{m-k}}
 \right\|_{L^{p}(\Omega,d^{\widetilde \alpha p})}
 \le
C \left\|
  {\nabla^ju}
 \right\|_{L^{p}(\Omega,d^{p-1})}.
\]
Altogether, inequality \eqref{april40} follows. The proof is
complete.
 \qed

\bigskip
\par\noindent
We conclude by demonstrating  the sharpness of Theorem
\ref{proposition:1.3}. Let us begin with condition \eqref{pq}.

\begin{proposition}\label{propo2}
Let $n, k,m\in \mathbb{N}$,  $n, m\ge 2$,
 and $1\le k\le m-1$.
 Assume  that
 $1\le p \le q <\infty$ and
\begin{equation}\label{pq>}
\frac 1q < \frac {n-p(m-k)}{np}.
\end{equation}
 Then inequality \eqref{eq:1.3} fails in any (smooth) bounded open set $\Omega \subset \mathbb R^n$.
\end{proposition}
\par\noindent
{\bf Proof}. We may assume, without loss of generality, that $0 \in
\Omega$. Suppose that inequality  \eqref{eq:1.3} holds for every $u
\in C^\infty _0(\Omega)$, with $m, k, p, q, r$ as in the statement.
Fix any such function $u$, and, for $\lambda
>1$,  consider the function $u_\lambda$ in
$\Omega$ defined as
\begin{equation}\label{lambda}
u_\lambda (x) = u(\lambda x) \quad \hbox{for $x \in \Omega$.}
\end{equation}
 Hence,
in particular,
\begin{equation}\label{april23}
\left( \int_{\Omega} d(x)^r\left|\nabla ^k\left(\frac{u_\lambda
(x)}{d(x)^{m-k}}\right)\right|^q\, dx \right)^{\frac{1}{q}}  \leq C
\sum_{j=0}^m\left( \int_{\Omega} d(x)^{p-1}\left|\nabla ^j{u_\lambda
(x)}\right|^p\, dx\right)^{\frac{1}{p}}
\end{equation}
for $\lambda >1$. Inequality \eqref{april23} can be rewritten, via a
change of variable, as
\begin{multline}\label{april24}
\lambda ^{k-\frac nq}\left( \int_{\Omega} d(x/\lambda)^r\left|\nabla
^k\left(\frac{u (x)}{d(x/\lambda)^{m-k}}\right)\right|^q\, dx
\right)^{\frac{1}{q}} \\ \leq C \sum_{j=0}^m \lambda ^{j-\frac
np}\left( \int_{\Omega} d(x/\lambda)^{p-1}\left|\nabla ^j{u
(x)}\right|^p\, dx\right)^{\frac{1}{p}}
\end{multline}
for $\lambda >1$. Since $\lim _{\lambda \to \infty }d(x/\lambda) =
d(0)$ and $\lim _{\lambda \to \infty }\nabla ^j
(d(x/\lambda)) = 0$ for $j = 1,
\dots , m$, uniformly in $x \in \Omega$, all the integrals in
\eqref{april24} converge to a finite limit as $\lambda \to \infty$.
Hence, passing to the limit  as $\lambda \to \infty$ in
\eqref{april24} leads to a contradiction, inasmuch as inequality
\eqref{pq>} is in force. \qed

\medskip
\par\noindent
The optimality of assumption \eqref{r} in Theorem
\ref{proposition:1.3} is the object of our last result.

\begin{proposition} \label{propo1}
 Let $n, k,m\in \mathbb{N}$,
$m\ge 2$,
 and $1\le k\le m-1$.
 Assume  that
 $1\le p \le q <\infty$, and
\begin{equation}\label{r<}
r< \frac qp (n-1) -n +q\,.
\end{equation}
Then there exists a smooth bounded open set $\Omega \subset \mathbb
R^n$ such that inequality \eqref{eq:1.3} fails.
\end{proposition}

%
\par\noindent
{\bf Proof}.
 Assume  first that $n\ge 2$.
Let $B^{n-1}_\rho (0)$ denote the ball in $\mathbb
R^{n-1}$, centered at $0$, with radius $\rho$. Let $\Omega$ be a
smooth bounded open set in $\mathbb R^n$ such that
$B_1^{n-1}(0)\times (0,1) \subset \Omega$, $B^{n-1}_\rho (0) \times
\{0\} \subset \partial \Omega$, and
$$d(x)=x_n \quad  \hbox{for $x\in B_1^{n-1}(0)\times (0,1)$.}$$
For instance, the set $\Omega$ can just be obtained by smoothing the
cylinder $B_2^{n-1}(0)\times (0,2)$.
Let $u$ be any smooth function in $\mathbb R^n$, compactly supported
in $B_1^{n-1}(0)\times (0,1)$. Hence, in particular, (the
restriction of) $u$ to $\Omega$ belongs to  $C_0^{\infty}(\Omega)$.
Given any $\lambda
>1$, consider the function $u_\lambda : \Omega \to \mathbb R$ defined as
 in \eqref{lambda}.
%
%
%
Since we are assuming that $\lambda >1$, we have that $u_\lambda \in
C_0^{\infty}(B_1^{n-1}(0)\times (0,1)) \subset
C_0^{\infty}(\Omega)$. Suppose that  inequality \eqref{eq:1.3}
holds for every $u \in C^\infty _0(\Omega)$, with  $m, k, p, q, r$
as in the statement. The choice of $u_\lambda$ as a trial function
in \eqref{eq:1.3} implies that
\begin{multline}\label{april21}
\left( \int_{B_1^{n-1}(0)\times (0,1)} x_n^{r}\left|\frac{\partial
^k}{\partial x_n^k}\left(\frac{u_\lambda
(x)}{x_n^{m-k}}\right)\right|^q\, dx \right)^{\frac{1}{q}} \\ \leq C
\sum_{j=0}^m\left( \int_{B_1^{n-1}(0)\times (0,1)}
x_n^{p-1}\left|\nabla ^j{u_\lambda (x)}\right|^p\,
dx\right)^{\frac{1}{p}}
\end{multline}
for $\lambda >1$. Hence, via a change of variable,
\begin{multline}\label{april22}
 \lambda^{-\frac rq
+m-\frac{n}{q}} \left( \int_{B_1^{n-1}(0)\times (0,1)}
x_n^{r}\left|\frac{\partial ^k}{\partial
x_n^k}\left(\frac{u(x)}{x_n^{m-k}}\right)\right|^q dx
\right)^{\frac{1}{q}}
\\ \leq C
\sum_{j=0}^m\lambda^{j-1+\frac 1p-\frac{n}{p}}\left(
\int_{B_1^{n-1}(0)\times (0,1)} {x_n}^{p-1}\left|\nabla
^j{u(x)}\right|^p dx\right)^{\frac{1}{p}}
\end{multline}
for $\lambda >1$. Letting $\lambda \to \infty$ in \eqref{april22}
yields a contradiction, under assumption \eqref{r<}.
\\
When $n=1$, the same argument applies on replacing
$B_1^{n-1}(0)\times (0,1)\subset \mathbb{R}^n$ with the interval
$(0,1)$.   \qed

\medskip
\par\noindent
{\bf Acknowledgments.} The authors are grateful to
Yoshinori Yamasaki for some helpful discussions.
\\ This research was initiated during a visit of the second named author
at the Department of Mathematics and Informatics ``U.Dini" of the
University of Florence, in the fall-winter semester 2013-2014. He
wishes to thank the members of the Department for their kind
hospitality.
\\ This work was partly funded by:  Research project of MIUR (Italian
Ministry of Education, University and Research) Prin 2012 ``Elliptic
and parabolic partial differential equations: geometric aspects,
related inequalities, and applications" (grant number 2012TC7588);
GNAMPA of the Italian INdAM (National Institute of High
Mathematics); JSPS KAKENHI (grant number 25220702).
%
%
%
%


\begin{thebibliography}{xx}


\bibitem[Ad1]{Ad1}
D.R. Adams,
Traces of potentials arising from translation invariant operators,
\emph{Ann. Sc. Norm. Super. Pisa} \textbf{25} (1971), 203--217.

\bibitem[Ad2]{Ad2}
D.R. Adams,
A trace inequality for generalized potentials, \emph{Studia Math.}
\textbf{48} (1973), 99--105.


\bibitem[CW]{CW}
H. Castro
\&
H. Wang, A Hardy type inequality for $W^{m,1}(0,1)$
functions,  \emph{Calc. Var.} \textbf{3-4} (2010), 525--531.

\bibitem[CDW1]{CDW1}
H. Castro, J. D{\'a}vila
\& H. Wang, 
A Hardy type inequality for $W^{2,1}_0(\Omega)$ functions,
\emph{C. R. Math. Acad. Sci. Paris} {\bf 349}
(2011), 765--767.


\bibitem[CDW2]{CDW2}
H. Castro, J. D\'avila \& H. Wang, A Hardy type inequality for
$W^{m,1}_0(\Omega)$ functions,  \emph{J. Eur. Math. Soc.} {\bf 15}
(2013), 145--155.



\bibitem[CEG]{CEG} A. Cianchi, D.E. Edmunds \& P.Gurka, On weighted
Poincar\'e inequalities, \emph{Math. Nachr.} \textbf{180} (1996),
15-41.


\bibitem[HKM]{HKM} J. Heinonen, T. Kilpel\"ainen \& O. Martio,  \emph{Nonlinear potential theory of degenerate elliptic
equations},  Oxford University Press, New York, 1993.


\bibitem[Ku]{K}
A. Kufner, \emph{Weighted Sobolev spaces}, Teubner-Texte zur
Mathematik, Leipzig, 1980.

\bibitem[Ho1]{H1}
T. Horiuchi,   The imbedding theorems for weighted Sobolev spaces,
\emph{J. Math. Kyoto Univ.} {\bf 29} (1989), 365--403.

\bibitem[Ho2]{H2}
T. Horiuchi,
 The imbedding theorems for weighted Sobolev spaces. II,
\emph{Bull. Fac. Sci. Ibaraki Univ. Ser. A} {\bf 23} (1991), 11--37.


\bibitem[Ma]{mazya}
V.G. Maz'ya, \emph{Sobolev spaces with applications to elliptic
partial differential equations},  Springer, Berlin, 2011.

\bibitem[MS]{MS}
M.K.V. Murthy \& G. Stampacchia,
 Boundary value problems for some degenerate-elliptic operators,
\emph{Ann. Mat. Pura Appl.} {\bf 80} (1968), 1--122.


\bibitem[OK]{OK} B. Opic \& A. Kufner,  \emph{Hardy-type inequalities},  Longman Scientific \& Technical, New York, 1990.



\end{thebibliography}
\end{document}